\documentclass{amsart}

\usepackage[margin=1.3in]{geometry}
\usepackage[utf8]{inputenc}
\usepackage[english]{babel}
\usepackage{csquotes}

\usepackage{amsmath}
\usepackage{amssymb}
\usepackage{physics}
\usepackage{youngtab}
\usepackage{mathtools}
\usepackage{comment}

\usepackage{graphicx} 
\usepackage{caption}
\usepackage{subcaption}
\usepackage{float}

\usepackage{algorithm,algpseudocode}
\usepackage{hyperref}
\usepackage{amsthm}
\usepackage[nameinlink]{cleveref}

\usepackage{nicematrix}
\usepackage{tikz}

\usepackage{todonotes} 
\usepackage[htt]{hyphenat} 
\usepackage{placeins}

\DeclareMathOperator{\pf}{Pf}
\DeclareMathOperator{\diag}{diag}
\DeclareMathOperator{\Ai}{Ai}

\newtheorem{proposition}{Proposition}[section]
\newtheorem{lemma}{Lemma}[section]
\theoremstyle{definition}
\newtheorem{definition}{Definition}[section]

\theoremstyle{remark}
\newtheorem{remark}{Remark}[section]

\makeatletter
\let\c@proposition\c@theorem
\let\c@lemma\c@theorem
\let\c@definition\c@theorem
\let\c@example\c@theorem
\let\c@xca\c@theorem
\let\c@remark\c@theorem
\let\c@corollary\c@theorem
\makeatother

\title[Sampling Pfaffian point processes]{Sampling  Pfaffian point processes 
and \\ the  symplectic Arnoldi method 
}

\author{Alan Edelman}
\address{Department of Mathematics, Massachusetts Institute of Technology}
\email{edelman@mit.edu}

\author{Sungwoo Jeong}
\address{Department of Mathematics, Cornell University}
\email{sjeong@cornell.edu}

\author{Simeon Schaub}
\address{Department of Physics, Karlsruhe Institute of Technology}
\email{simeon@schaub.rocks}

\date{}

\begin{document}

\begin{abstract}
We present an exact sampling algorithm for Pfaffian point processes based on a skew-symmetric analogue of the Cholesky factorization. This algorithm enables efficient sampling of a wide range of statistics arising in random matrix theory and combinatorics.
For instance, we can sample eigenvalues of the orthogonal and symplectic ensembles ($\beta = 1,4$).

In addition, we introduce a symplectic Arnoldi method for computing skew-orthogonal polynomials associated with a general weight function. This method can be used to efficiently construct the $2 \times 2$ matrix valued skew-symmetric kernels that arise in $\beta = 1,4$ polynomial ensembles. We illustrate our approach with several numerical examples and experiments, including the symmetric corner growth model, the finite-$N$ Gaussian (Hermite) orthogonal and symplectic ensembles, and the $\beta = 1,4$ Airy point processes and Tracy–Widom distributions.
\end{abstract}

\maketitle

\section{Introduction}

Pfaffian point processes (PfPPs) arise in many areas of mathematics, including random matrix theory, statistics, and combinatorics. We say that a point process is a Pfaffian point process if its $k$-point correlation function $\rho_k$ satisfies 
\begin{equation*}
    \rho_k(x_1, \dots, x_k) = \pf[K(x_i, x_j)]_{i, j \in [k]},
\end{equation*}
for some $2\times 2$ matrix valued kernel $K$. Most famously, the eigenvalues of orthogonal and symplectic ensembles ($\beta=1,4$) have Pfaffian correlation functions \cite{mehta2004random}. Several other stochastic processes, such as random involutions \cite{rains2000correlation,forrester2006correlation,baik2018pfaffian}, symmetric corner growth \cite{johansson_shape_2000}, one-dimensional particle systems \cite{garrod2018examples,sniady2026determinant}, or zeros of a Gaussian power series \cite{matsumoto2013correlation}, also have correlation functions expressed as Pfaffians.

A closely related and more widely studied point process is the determinantal point process (DPP). A point process is determinantal if its $k$-point correlation function has a determinantal expression. Determinantal point processes describe many interesting combinatorial objects and random matrix statistics, such as the eigenvalues of orthogonal polynomial ensembles \cite{mehta2004random}, dimer models \cite{kasteleyn1963dimer,kenyon1997local}, non-intersecting paths \cite{johansson2002non}, uniform measure of spanning trees \cite{burton1993local}, and the $\beta = 2$ Airy point process and Tracy-Widom distribution \cite{tracy1994level}. It has also found applications in numerical linear algebra \cite{derezinski2021determinantal} and machine learning \cite{kulesza2012determinantal}. 

For DPPs, dozens of sampling algorithms have been introduced in the last two decades \cite{li2016fast,anari2016monte,derezinski2019exact,bardenet2020monte}. In particular, in 2006, Hough et al.\ introduced one of the first sampling algorithms for DPPs \cite{hough2006determinantal}. This algorithm is exact and uses the eigendecomposition of the $K$ kernel (or the $L$ kernel), which restricts the algorithm to symmetric DPPs. More recently, another exact DPP sampling algorithm was introduced \cite{poulson2020high,launay2020exact}. This second algorithm uses conditional probabilities and the LU factorization, and it enables sampling from general (nonsymmetric) DPPs, such as the Aztec diamond domino tiling \cite{jockusch1998random,johansson2002non,chhita2015asymptotic} and the Airy process \cite{prahofer2002scale,edelman2023conditional}.

On the other hand, for PfPPs, only a limited number of sampling algorithms have been developed \cite{Bardenet_2024}. In this work, we introduce a new sampling algorithm for PfPPs, which is exact. Our algorithm uses conditional measures of Pfaffian point processes \cite{borodin2005eynard,bufetov2021conditional}, together with a skew-symmetric variant of the Cholesky factorization. Our algorithm can be thought of as an analogue of the second DPP sampling algorithm mentioned above \cite{launay2020exact,poulson2020high} and it can sample from any PfPP. 

Note that a PfPP sampling algorithm can be used to sample various random matrix eigenvalues or random measures arising in combinatorics. As previously mentioned, one can use a PfPP sampling algorithm to sample from a rank $N$ $\beta=1,4$ polynomial ensemble, i.e.\ a point process $X$ with fixed $|X|=N$ and the joint probability distribution for $X = \{x_1, \dots, x_N\}$,
\begin{equation}\label{eq:betaensemble}
    f(X) = \frac{1}{\mathcal{Z}}\prod_{i<j} |x_i - x_j|^\beta \prod_{j=1}^N w(x_j),
\end{equation}
where $w(x)$ is the given weight function and $\mathcal{Z}$ is the appropriate partition function (normalization). The density $f(X)$ can be expressed as a Pfaffian of a rank-$N$ kernel with $2\times 2$ matrix entries. For an example see \cite{tracy1996orthogonal,mehta2004random,a2005matrix}. More precisely, one needs to compute a set of skew-orthogonal polynomials with respect to a given weight $w(x)$ to obtain the kernel. Recall that the Lanczos (a special case of Arnoldi) iteration can be used to compute orthogonal polynomials through their three-term recurrence \cite{golub1969calculation,gautschi2004orthogonal}. To efficiently compute the skew-orthogonal polynomials for $\beta=1,4$ Pfaffian point process kernels, we propose a generalization of the Arnoldi iteration, which can be used to compute the skew-orthogonal polynomials and thus the symmetric kernel for PfPPs. 

To demonstrate the effectiveness of our algorithms, we provide various numerical experiments. First, we sample from the symmetric corner growth process, which is a known PfPP \cite{johansson_shape_2000}. Then, we compute the eigenvalues of the finite-$N$ GOE and GSE, and compare our algorithm with other sampling algorithms. Since many applications are continuous PfPPs, we discuss various methods to handle the discretization, or avoid discretization altogether. Finally, we perform numerical experiments on the soft-edge scaling limits of the GOE and GSE, also as known as the Tracy-Widom distribution with $\beta=1,4$, or the Airy point processes. Furthermore by using the conditional Pfaffian point process approach \cite{bufetov2021conditional,edelman2023conditional}, we sample various statistics of the $k$ largest eigenvalues of the Airy point process. 

This work is organized as follows: In \Cref{sec:pfppsampler}, we briefly review the basics of the Pfaffian and PfPPs, and introduce an exact sampling algorithm for PfPPs. We illustrate examples and applications of PfPPs in random matrix theory and combinatorics. In \Cref{sec:skewarnoldi}, we propose a symplectic Arnoldi iteration, that can be used to efficiently compute the skew-orthogonal polynomials associated to $\beta=1,4$ polynomial ensembles and corresponding skew-symmetric kernels. Finally, in \Cref{sec:numerical}, we illustrate various numerical experiments.

\section{Pfaffian point process and its sampling algorithm}\label{sec:pfppsampler}

\subsection{Pfaffian Point Process}

The Pfaffian is a kind of analogue of the determinant for 
real or complex $2n\times 2n$ skew-symmetric matrices, $A = -A^T$.
In fact it is $\pm  \sqrt{\det A},$ but the tricky part is 
which square root to take. 
The precise definition that picks a sign is
as follows:
\begin{equation}
    \pf(A) \coloneq \frac{1}{2^n n!}\sum_{\sigma\in S_{2n}} \text{sign}(\sigma)\prod_{i=1}^n A_{\sigma(2i-1), \sigma(2i)},
\end{equation}
where $S_{2n}$ is the set of all permutations of $[2n]=\{1, \dots, 2n\}$. 

Let $K$ be a $2n\times 2n$ matrix and $X, Y\subseteq [n]$ be two subsets of $[n]$ with sizes $|X|=j, |Y|=k$. We denote by $K_{X,Y}$ the $2j \times 2k$ submatrix of $K$ obtained by collecting rows $\{2x-1 \,|\, x\in X\} \cup \{2x \,|\, x \in X\}$ and columns $\{2y-1\, | \,y \in Y\}\cup\{2y \,|\, y\in Y\}$. This indexing is natural if $K$ is viewed as a $n\times n$ matrix with each entry being a $2\times 2$ matrix. We will simply write $K_X$ for $K_{X, X}$. 
The following $2n\times 2n$ matrix $J_n$ plays the role of the identity matrix for $2n\times 2n$ skew-symmetric matrices
\begin{equation}\label{eq:Jndef}
    J_n \coloneq I_n \otimes \begin{bNiceArray}{rc}
    0 & 1 \\ -1 & 0
    \end{bNiceArray} = \begin{bNiceArray}{rccrc}
        0 & 1 &  & & \\
        -1 & 0 &  & & \\
        & & \ddots & & \\
        & & & 0 & 1\\
        & & & -1 & 0
    \end{bNiceArray}.
\end{equation}

We define the \textit{point process} to be a probability measure on all (locally finite) collections of points $X$ selected from some set $\mathcal{G}$. The set $\mathcal{G}$ is sometimes called the ground set. A (discrete) Pfaffian point process can be defined as follows:

\begin{definition}[Discrete Pfaffian point process]
Let $\mathcal{G}$ be a discrete set. A point process $\mathcal{J}$ is a (discrete) Pfaffian point process if there exists a $2|\mathcal{G}|\times 2|\mathcal{G}|$ skew-symmetric matrix $K$, such that
\begin{equation}
    \mathbb{P}(S\subset \mathcal{J}) = \pf(K_S),
\end{equation}
for any subset $S$ of $\mathcal{G}$.
\end{definition}
The matrix $K$ is often regarded as a ($2\times 2$ matrix valued) matrix indexed by $\mathcal{G}$. The matrix $K$ is called the $K$-\textit{kernel} or simply the \textit{kernel}, to distinguish from the $L$-\textit{kernel}, which is introduced below. Note that a PfPP may not have an $L$-kernel. 

\begin{definition}[$L$-kernel]
The $L$-kernel of a PfPP $\mathcal{J}$ is defined as a $2n\times 2n$ matrix $L$ that satisfies
\begin{equation*}
    \mathbb{P}(\mathcal{J} = S) = \frac{\pf(L_S)}{\pf(J_n+L)},
\end{equation*}
for all subsets $S$ of $[n]$.
\end{definition}
For DPPs, there is a well-known relationship between $L$ and $K$, namely $L = K(I-K)^{-1}, K = L(I+L)^{-1}$. For PfPPs, we have similar relationships \cite{borodin2005eynard},
\begin{gather*}
    L = J_nK(J_n - K)^{-1} = (K-J_n)^{-1} - J_n,\\ K = J_n + (J_n + L)^{-1} = J_n L (J_n+L)^{-1}.
\end{gather*}

To define continuous Pfaffian point processes, we let $\mathcal{G}$ be a continuous set. First we define the $k$-point correlation function $\rho_k$ as a function of $k$ variables, such that
\begin{equation*}
    \mathbb{E}\left[\sum_{x_1, \dots, x_k\in\mathcal{J}}f(x_1, \dots, x_k)\right] = \int_{\mathcal{G}^k} f(x_1, \dots, x_k)\rho_k(x_1, \dots, x_k)\dd{x_1}\cdots \dd{x_k}
\end{equation*}
holds for any nonnegative measurable $f:\mathcal{G}^k\to \mathbb{R}^+$. Equivalently, the $k$-point correlation function can also be defined as 
\begin{equation*}
    \rho_k(x_1, \dots, x_k) = \lim_{\dd{x}\to 0}\frac{\mathbb{P}\left(\text{$k$ points in $[x_1, x_1+\dd{x}], \dots, [x_k, x_k+\dd{x}]$}\right)}{(\dd{x})^n}. 
\end{equation*}
See related discussions in \cite{tracy1998correlation,mehta2004random} for further details.

A continuous point process $\mathcal{J}$ on $\mathcal{G}$ is Pfaffian when its $k$-point correlation function can be expressed as a Pfaffian of $2\times 2$ matrix valued skew-symmetric kernel $K:\mathcal{G}\times \mathcal{G} \to \mathbb{C}^{2\times 2}$.
\begin{definition}[Continuous Pfaffian point process]
A point process $\mathcal{J}$ is a (continuous) Pfaffian point process if its $k$-point correlation function $\rho_k$ is given as
\begin{equation}
    \rho_k(x_1, \dots, x_k) = \pf\left[K(x_i, x_j)\right]_{i, j\in [k]},
\end{equation}
where $K:\mathcal{G}\times \mathcal{G}\to\mathbb{C}^{2\times 2}$ is a $2\times 2$ matrix valued kernel and $\left[K(x_i, x_j)\right]_{i, j\in [k]}$ is viewed as a $2k\times 2k$ skew-symmetric matrix.
\end{definition}

Pfaffian point processes are sometimes also called Pfaffian random fields. For more details on PfPPs, we refer to \cite{rains2000correlation,soshnikov2003janossy,borodin2005eynard}.

\subsubsection{Random Matrix Theory}

Consider a weight function $w(x)$ and associated orthogonal polynomials $\{P_n(x)\}_{n\in\mathbb{N}}$ with respect to the inner product $\langle f, g, \rangle_w = \int f(x)g(x)w(x)\dd{x}$. The joint probability density of the orthogonal polynomial ensemble associated with $w(x)$ is given by \eqref{eq:betaensemble} with $\beta=2$. It can be shown that the joint density equals 
\begin{equation*}
    f(x_1, \dots, x_N) = \det\left[K(x_i, x_j)\right]_{i, j\in[N]},
\end{equation*}
where the kernel $K$ is given as
\begin{equation*}
    K_N(x, y) \coloneq \sqrt{w(x)w(y)}\sum_{j=0}^{N-1} \frac{P_j(x)P_j(y)}{h_j},
\end{equation*}
where $h_j = \langle P_j, P_j \rangle_w$. Using the Christoffel-Darboux formula, this simplifies to 
\begin{equation*}
    K_N(x, y) = \frac{\sqrt{w(x)w(y)}}{h_{n-1}}\cdot \frac{\kappa_{n-1}}{\kappa_n} \frac{P_N(x)P_{N-1}(y) - P_N(y) P_{N-1}(x)}{x-y},
\end{equation*}
where $\kappa_n$ is the leading coefficient of $P_n$. Gaudin and Mehta \cite{mehta1960density,gaudin1961loi,mehta2004random} noticed that for orthogonal polynomial ensembles, the joint density and $N$-point correlation functions are equivalent, as it turns out that this kernel satisfies $K^2 = K$ in the sense of the integral operator. Thus, the kernel $K_N$ defines a determinantal point process. 

For $\beta=1,4$, we can obtain similar results using Pfaffians instead of determinants. Consider the rank $N$ $\beta=1$ or $\beta=4$ polynomial ensemble corresponding to $w(x)$, defined by \eqref{eq:betaensemble}. In both cases, one can use the skew-orthogonal polynomials associated with $w$ to obtain 
\begin{equation*}
    f(x_1, \dots, x_N) = \pf\left[K(x_i, x_j)\right]_{i, j\in[N]},
\end{equation*}
where $K$ is a $2\times 2$ matrix valued kernel, which we will define in \Cref{sec:GOEfinite,sec:GSEfinite}. For general $\beta>0$, it is an open question whether the eigenvalues of G$\beta$E can be explained by special point processes similar to determinantal or Pfaffian point processes. 

Throughout this work, we frequently use the Gaussian (Hermite) ensemble as a representative example. For any $\beta>0$, the G$\beta$E is defined as a system of $N$ particles with joint probability density
\begin{equation}\label{eq:hermitebeta}
    f(\lambda_1, \dots, \lambda_N) = \frac{1}{\mathcal{Z_{N, \beta}}} \prod_{i<j}|x_i - x_j|^\beta \prod_{i=1}^N e^{-\frac{\beta}{4}\lambda_i},
\end{equation}
where 
\begin{equation}\label{eq:partitionfunction}
    \mathcal{Z_{N, \beta}} = (2\pi)^{\frac{N}{2}}(2^{-1}\beta)^{-\frac{N}{2}-\frac{\beta}{4}N(N-1)}\prod_{j=1}^N \frac{\Gamma(1 + \frac{\beta j}{2})}{\Gamma(1 + \frac{\beta}{2})}.
\end{equation}
For $\beta=1,2,4$, the joint density of the eigenvalues of random Gaussian symmetric, Hermitian and self-dual matrices are equal to \eqref{eq:hermitebeta}. For general $\beta>0$, the joint density \eqref{eq:hermitebeta} can be obtained from the eigenvalues of certain tridiagonal matrix models \cite{dumitriu2002matrix}.

\subsection{Sampling algorithm for PfPPs using the Cholesky-like decomposition}

In this section, we propose a sampling algorithm for discrete PfPPs. Our algorithm, \Cref{alg:pfppsampler}, uses the Cholesky-like decomposition for skew symmetric matrices \cite{benner2000cholesky}.

\begin{lemma}[Cholesky-like decomposition for skew-symmetric matrices]\label{lem:bjb}
    For a given real $2n\times 2n$ skew-symmetric matrix $A$, there exists a matrix $B\in\mathbb{R}^{2n\times 2n}$ such that $A = BJ_n B^T$ holds. 
    Similarly, for any complex $2n\times 2n$ skew-symmetric matrix $A$, there exists a matrix $B\in\mathbb{C}^{2n\times 2n}$ such that $A = BJ_nB^T$.
\end{lemma}

The following conditional probabilities \cite{borodin2005eynard} related to PfPPs are useful:
\begin{proposition}\label{prop:condppp}
    Let $K\in\mathbb{C}^{2n\times 2n}$ be the kernel of a discrete PfPP $\mathcal{J}$ on the ground set $\mathcal{G}$. Given disjoint subsets $X, Y\subset \mathcal{G}$, the following conditional probabilities have Pfaffian expressions:
    \begin{align}\label{eq:cond1}
        \mathbb{P}(X \subset \mathcal{J}\,|\, Y \subset \mathcal{J}) &= \pf\left(K_{X}-K_{X,Y}K_{Y}^{-1}K_{Y,X}\right), \\ \label{eq:cond2}
        \mathbb{P}(X \subset \mathcal{J}\,|\, Y \subset \mathcal{J}^c) &= \pf\left(K_{X}-K_{X,Y}(K_{Y} - J_{|Y|})^{-1}K_{Y,X}\right),
    \end{align}
    where $J_{|Y|}\in \mathbb{R}^{2|Y|\times 2|Y|}$ is defined as in \eqref{eq:Jndef}. Let $Y' = \mathcal{G}\backslash Y$. 
\end{proposition}

If we define two matrices $K^{\text{in}}, K^{\text{out}}$ of order $2n-2|Y|$,
\begin{equation*}
    K^{\text{in}} = K_{Y'} - K_{Y',Y}K_Y^{-1}K_{Y,Y'}, \hspace{0.5cm}
    K^{\text{out}} = K_{Y'} - K_{Y',Y}(K_Y-J_{|Y|})^{-1}K_{Y,Y'},
\end{equation*}
we can simply write \eqref{eq:cond1}, \eqref{eq:cond2} as
\begin{equation*}
    \mathbb{P}(X \subset \mathcal{J}\,|\, Y \subset \mathcal{J}) = \pf(K^{\text{in}}_{X}), \hspace{1cm} \mathbb{P}(X \subset \mathcal{J}\,|\, Y \subset \mathcal{J}^c) = \pf(K^{\text{out}}_{X}).
\end{equation*}

Since we can arbitrarily select $X\subset Y'$, \Cref{prop:condppp} shows that the point process $\mathcal{J}$, restricted to $Y'$ with the condition $Y\subset\mathcal{J}$, is a PfPP with kernel $K^{\text{in}}$. See \cite{bufetov2021conditional} for a more rigorous derivation. Similarly, the point process $\mathcal{J}$, restricted to $Y'$ and conditioned with $Y\subset\mathcal{J}^c$, is a PfPP with kernel $K^{\text{out}}$. In particular, when $Y$ is a set containing a single particle $\{i\}$, we have $\mathbb{P}(X\subset \mathcal{J}\,|\,Y\subset J^c) = \mathbb{P}(X\subset\mathcal{J}\,|\,i\notin \mathcal{J})$.

Our sampling algorithm follows from the above Proposition. For simplicity, assume $\mathcal{G} = [n]$ and let $O=\{1\}$, $O'=\{2, \ldots, n\}$. Let us begin with the first index (particle) $1$. Since the probability $\mathbb{P}(1\in\mathcal{J}) = \pf(K_{O})$ equals the entry at index $(1,2)$ of $K$ (viewed as $2n\times 2n$ matrix), we perform a Bernoulli trial with success probability equal to this entry to decide whether to include $1$ in the sample or not. If $1$ is determined to be included in the sample (accepted), we derive a new kernel from \Cref{prop:condppp}:
\begin{equation*}
    K^\text{in} = K_{O'}- K_{O',O}K_{O}^{-1}K_{O,O'}.
\end{equation*}
Then the PfPP with kernel $K^{\text{in}}$ specifies the distribution of the rest of the sample. Therefore, we will return $\{1\} \cup \mathcal{J}^{\text{in}}$, where $\mathcal{J}^\text{in}$ is a PfPP with $K^\text{in}$ and ground set $\{2, \dots, n\}$. 

When $1$ is rejected, we similarly define
\begin{equation*}
    K^\text{out} = K_{O'}- K_{O'O}(K_{O}-J_1)^{-1}K_{OO'}.
\end{equation*}
We can proceed with the PfPP $\mathcal{J}^\text{out}$ with kernel $K^\text{out}$. Iteratively applying this step to $O = \{2\}, \{3\}, \dots,$ we obtain a sample. \Cref{alg:pfppsampler} can be thought of as a generalization of the DPP sampling algorithm introduced in \cite{launay2020exact,poulson2020high}.

\begin{algorithm}[h]
\begin{algorithmic}[1]
\Function{randPfPPcholesky}{$K \in \mathbb{R}^{2n \times 2n}$}
    \State $\mathcal{I} \gets [\,]$
    \For{$j = 1, \dots, n$}
        \State $p \gets K_{2j - 1, 2j}$ \Comment{Pfaffian of $K_{2j - 1:2j,\ 2j - 1:2j}$}
        \State $b \sim \operatorname{Bernoulli}(p)$
        \If{b}
            \State $\mathcal{I} \gets [\mathcal{I} \,|\, j]$
        \Else
            \State $p \gets p - 1$
        \EndIf
        \State $K_{2j + 1:2n,\ 2j + 1:2n} \gets K_{2j + 1:2n,\ 2j + 1:2n} - K_{2j + 1:2n,\ 2j - 1:2j} \cdot p^{-1} J \cdot K_{2j - 1:2j,\ 2j + 1:2n}$
\[
\scalebox{.7}{$\begin{bNiceMatrix}
\phantom{1} & \phantom{1} & \phantom{1} & \phantom{1} & \phantom{1} & \phantom{1} \\
\\ \\ \\ \\
\phantom{1}
\CodeAfter
\tikz \draw [style=dotted,line width=1pt] (1.5-|1.5) -- (2.7-|2.7) -- (2.7-|1-last) -- (1.5-|1-last) -- cycle ;
\tikz \draw [fill=red!50] (4-|4) -- (last-1-|4) -- (last-1-|1-last) -- (4-|1-last) -- cycle ;
\tikz{
  \draw[line cap=round,line width=2pt,
        dash pattern=on 0pt off 9pt]
    (3.3-|3.7) -- (3.3-|1-last)
    (3.7-|4.3) -- (3.7-|1-last);
  \draw[line cap=round,line width=2pt,
        dash pattern=on 0pt off 7.2pt]
    (3.7-|3.3) -- (last-1-|3.3)
    (4.3-|3.7) -- (last-1-|3.7) ;
}
\end{bNiceMatrix}$}
\gets
\scalebox{.7}{$\begin{bNiceMatrix}
\phantom{1} & \phantom{1} & \phantom{1} & \phantom{1} & \phantom{1} & \phantom{1} \\
\\ \\ \\ \\
\phantom{1}
\CodeAfter
\tikz \draw [style=dotted,line width=1pt] (1.5-|1.5) -- (2.7-|2.7) -- (2.7-|1-last) -- (1.5-|1-last) -- cycle ;
\tikz \draw [fill=red!50] (4-|4) -- (last-1-|4) -- (last-1-|1-last) -- (4-|1-last) -- cycle ;
\tikz{
  \draw[line cap=round,line width=2pt,
        dash pattern=on 0pt off 9pt]
    (3.3-|3.7) -- (3.3-|1-last)
    (3.7-|4.3) -- (3.7-|1-last);
  \draw[line cap=round,line width=2pt,
        dash pattern=on 0pt off 7.2pt]
    (3.7-|3.3) -- (last-1-|3.3)
    (4.3-|3.7) -- (last-1-|3.7) ;
}
\end{bNiceMatrix}$}
-
\scalebox{.7}{$\begin{bNiceMatrix}
\phantom{1} & \phantom{1} & \phantom{1} & \phantom{1} & \phantom{1} & \phantom{1} \\
\\ \\ \\ \\
\phantom{1}
\CodeAfter
\tikz \draw [style=dotted,line width=1pt] (1.5-|1.5) -- (2.7-|2.7) -- (2.7-|1-last) -- (1.5-|1-last) -- cycle ;
\tikz \draw [fill=green!50] (4-|3) -- (last-1-|3) -- (last-1-|4) -- (4-|4) -- cycle ;
\tikz{
  \draw[line cap=round,line width=2pt,
        dash pattern=on 0pt off 9pt]
    (3.3-|3.7) -- (3.3-|1-last)
    (3.7-|4.3) -- (3.7-|1-last);
  \draw[line cap=round,line width=2pt,
        dash pattern=on 0pt off 7.2pt]
    (3.7-|3.3) -- (last-1-|3.3)
    (4.3-|3.7) -- (last-1-|3.7) ;
}
\end{bNiceMatrix}$}
\cdot
\begin{bNiceMatrix}
& p^{-1}\\ -p^{-1}
\end{bNiceMatrix}
\cdot
\scalebox{.7}{$\begin{bNiceMatrix}
\phantom{1} & \phantom{1} & \phantom{1} & \phantom{1} & \phantom{1} & \phantom{1} \\
\\ \\ \\ \\
\phantom{1}
\CodeAfter
\tikz \draw [style=dotted,line width=1pt] (1.5-|1.5) -- (2.7-|2.7) -- (2.7-|1-last) -- (1.5-|1-last) -- cycle ;
\tikz \draw [fill=blue!50] (3-|4) -- (4-|4) -- (4-|1-last) -- (3-|1-last) -- cycle ;
\tikz{
  \draw[line cap=round,line width=2pt,
        dash pattern=on 0pt off 9pt]
    (3.3-|3.7) -- (3.3-|1-last)
    (3.7-|4.3) -- (3.7-|1-last);
  \draw[line cap=round,line width=2pt,
        dash pattern=on 0pt off 7.2pt]
    (3.7-|3.3) -- (last-1-|3.3)
    (4.3-|3.7) -- (last-1-|3.7) ;
}
\end{bNiceMatrix}$}
\]
    \EndFor
    \State \Return $\mathcal I$
\EndFunction
\end{algorithmic}
\caption{Pfaffian Point Process Sampler}\label{alg:pfppsampler}
\end{algorithm}

\subsection{Sampling from Continuous PfPPs}\label{sec:continuous}

Recall that many applications of the PfPP are continuous. For instance, eigenvalues of the GOE form a subset of $\mathbb{R}$, which is a continuous ground set. In this section, we will discuss strategies for sampling continuous PfPPs. 

The easiest way to sample from a continuous PfPP is to simply discretize the kernel in advance. One chooses a grid $X \times X$ with $X = \{x_\text{min}, x_\text{min} + \Delta, \dots, x_\text{max} - \Delta, x_\text{max}\}$ and evaluates $K$ at each point on this grid and multiplies the resulting matrix by $\Delta$. One can then sample from this matrix just like from a discrete PfPP and convert the sampled index set back to points on the grid.

The difficulty of this method is the choice of grid. It has to be extensive enough with small enough $\Delta$ to faithfully represent the continuous PfPP, but more and finer grid points make constructing the kernel matrix and sampling from it increasingly expensive, so a tradeoff has to be made. The main disadvantage however is that one can only ever sample points on the original grid, so the continuous nature of a PfPP is not captured well.

The rest of the methods discussed assume $K$ describes a rank-$N$ projection PfPP, so each sample has a known fixed size.

One way to avoid discretizing the whole kernel is the \textsc{RPCholesky} method introduced by Epperly and Moreno in \cite{epperly2023kernel}, which we modified for use with Pfaffians to use the skew-Cholesky factorization. The two difficulties with that method were finding an efficient way to sample from the diagonal -- e.g.\ due to the fact that $K(x, x) / N$ is not log-concave for the GSE, adaptive rejection sampling could not be used there -- and the fact we need a good upper bound on the conditional probabilities for the rejection sampling step.

It is also possible to employ more general Markov chain Monte Carlo (MCMC) methods to generate representative samples from a PfPP. Multiple such methods were studied as part of this work and their effectiveness verified on the GSE:

\begin{itemize}
    \item \textsc{Slice-within-Gibbs} \cite{neal1997slice}: One starts with an initial set of $N$ points. These could be sampled from a discretized version of the PfPP, but in our experiments, starting from equispaced points inside the support also worked well. At each step, one goes through all points $s_i$ and for each, samples a new value conditioned on all the other points $S_{-i}$ from $\pi(s_i | S_{-i})$ using slice sampling. The conditional probabilities were derived above. Similar to \mbox{\textsc{RPCholesky}}, one can work with the lower triangular skew-Cholesky factor $L$ of $K(S_{-i}, S_{-i})$, which is also cheaper to invert.
    \item \textsc{MALA-within-Gibbs} \cite{tong2020mala} works just as \textsc{Slice-within-Gibbs}, but the slice sampling step is replaced with a proposed update that is either accepted or rejected using the Metropolis-adjusted Langevin algorithm (MALA). The additional difficulty comes from requiring gradients for $\pi(s_i | S_{-i})$, which can be calculated using automatic differentiation. In this case, the \texttt{ForwardDiff.jl} \cite{RevelsLubinPapamarkou2016} library was used through \texttt{DifferentiationInterface.jl} \cite{dalle2025commoninterfaceautomaticdifferentiation,schafer2022abstractdifferentiationjlbackendagnosticdifferentiableprogramming,dalle_2026_18445764}. Steps were cheaper, but variance was higher for the same number of steps due to the fact that at each step, not all samples always get updated. This means \textsc{MALA-within-Gibbs} is slower in exploring the sample space, so, at least for the GSE specifically, \textsc{Slice-within-Gibbs} should be preferred.
    \item \textsc{Hamiltonian Monte Carlo} using the \textsc{NUTS} sampler \cite{hoffman2014no} was also studied. In contrast with the previous two methods, the joint PDF was calculated directly from the Pfaffian of $K(S, S)$. \texttt{Mooncake.jl} \cite{Tebbutt_Mooncake_Towards_a}, again through \texttt{DifferentiationInterface.jl}, was used for calculating the gradients. A custom rule for the log-Pfaffian was required, which we hope to upstream. For the overall sampling \texttt{DynamicHMC.jl} \cite{tamas_k_papp_2026_18130162} was used. This was the most expensive among the methods studied.
\end{itemize}

\section{Skew-orthogonal polynomial ensembles}

To construct PfPPs from a discrete $\beta = 1, 4$ Coulomb gas \eqref{eq:betaensemble}, one first needs to construct skew-orthogonal polynomials from the weight function $w(x)$. For $\beta = 1$, define a skew-symmetric inner product $\langle f, g \rangle_w^{(1)}$ as
\begin{equation}\label{eq:beta1innerproduct}
    \langle f, g \rangle_w^{(1)} \coloneq \sum_{x \in \mathcal{D}, y \in \mathcal{D}} f(x) g(y) \cdot \varepsilon(x, y) w(x) w(y), 
\end{equation}
where $\varepsilon(x, y) \coloneq \frac{1}{2} \operatorname{sign}(y - x)$. Similarly, for $\beta = 4$, define another inner product \cite{adler_classical_nodate}
\begin{equation}\label{eq:beta4innerproduct}
    \langle f, g \rangle_w^{(4)} \coloneq \sum_{x \in \mathcal{D}} \left( f(x) g'(x) - f'(x) g(x) \right) \cdot w(x)^2. 
\end{equation}
The goal now is to find (skew-orthogonal) polynomials $p_n$ of degree $n$, such that for all $m, n \in \mathbb{N}$
\begin{align*}
    \langle p_{2m}, p_{2n + 1} \rangle = -\langle p_{2n + 1}, p_{2m} \rangle &= \delta_{m, n},\\
    \langle p_{2m}, p_{2n} \rangle = \langle p_{2m + 1}, p_{2n + 1} \rangle &= 0,
\end{align*}
where the inner product is chosen from one of the two inner products defined above. We denote the skew-orthogonal polynomials associated with \eqref{eq:beta1innerproduct} by $\{R_n\}_{n\in\mathbb{N}}$, and the skew-orthogonal polynomials associated with \eqref{eq:beta4innerproduct} by $\{Q_n\}_{n\in\mathbb{N}}$, following the notation of \cite{mehta2004random}.

While these skew-orthogonal polynomials can be derived by hand from their orthogonal counterparts in some cases \cite{felipe_skew-orthogonal_2006, adler_classical_nodate}, there is a need for a more general numerical method that works for any weight function $w(x)$ (or, in fact, even arbitrary skew-symmetric inner products).

\subsection{Symplectic Arnoldi iteration}\label{sec:skewarnoldi}

Symplectic Arnoldi iteration has previously been studied in the context of Hamiltonian systems and ODEs \cite{fassbender_2007,celledoni2016symplectic}, focusing on the special case of symplectic Lanczos. Skew-orthogonalization/symplectic Gram-Schmidt goes back even further in the form of an SR factorization and numerical convergence has been studied by Watkins and Elsner \cite{watkins1991convergence}. To the best of our knowledge, the generalization of symplectic Lanczos to a more general Arnoldi iteration has not been studied extensively before.

Salam \cite{salam2005theoretical} does already talk about an extension of the Lanczos method for orthogonal polynomials to the skew-orthogonal, or what they call the symplectic, case:
\begin{quote}
    Other investigations could be pursued, as, for example, introducing symplectic polynomials in contrast with orthogonal polynomials for Lanczos methods.
\end{quote}
However, it seems like this may not been pursued any further. Therefore, inspired by the Lanczos method for calculating coefficients of the three-term recurrence for orthogonal polynomial ensembles \cite{gautschi2004orthogonal,qu_lanczos_2025}, we present an analogous scheme for skew-orthogonal polynomials. The trick for the Lanczos method is to tridiagonalize the operator
\begin{equation*}
    A: f(x) \longmapsto x \cdot f(x),
\end{equation*}
with respect to the symmetric inner product
\begin{equation*}
    \langle f, g \rangle_w^2 \coloneq \sum_{x \in \mathcal{D}} f(x) g(x) \cdot w(x).
\end{equation*}
By iteratively applying the map $A$ to an initial constant function and employing Gram-Schmidt orthogonalization to orthogonalize the result with respect to all previous functions.

While it is easy to see that the operator $A$ is Hermitian with respect to the symmetric inner product $\langle \cdot, \cdot \rangle_w^2$, e.g. $\langle f, A g \rangle_w^2 = \langle A f, g \rangle_w^2$, this does not hold anymore for the previously mentioned skew-symmetric inner products in \eqref{eq:beta1innerproduct} and \eqref{eq:beta4innerproduct}. This also implies that skew-orthogonal polynomials no longer fulfill a three-term recurrence
\begin{equation*}
    x \cdot p_n(x) = \beta_{n + 1} p_{n + 1}(x) + \alpha_n p_n(x) + \beta_n p_{n - 1}(x).
\end{equation*}
Instead, $x \cdot p_n(x)$ is a linear combination of all polynomials up to degree $n + 1$ in general:
\begin{equation}\label{eq:skew_recurrence}
    x \cdot p_n(x) = \sum_{k = 0}^{n + 1} H_{k, n} p_k(x).
\end{equation}

Since the map $A$ is not any more Hermitian, what was previously the Lanczos iteration now becomes the more general Arnoldi iteration, but we now also require our basis to be symplectic instead of orthogonal, meaning at the $M^\text{th}$ iteration, we are looking for a symplectic basis $S \in V^M$ of our function space $V$, an upper Hessenberg matrix $H \in \mathbb{R}^{M \times M}$ and a residual $r \in V$, such that
\begin{equation*}
    A S_n = \sum_{k = 0}^{\min(n + 1, M - 1)} S_k H_{k, n} + \delta_{n, M - 1} r,
\end{equation*}
with the $S_n$ fulfilling
\begin{align*}
    \langle S_{2m}, S_{2n + 1} \rangle = -\langle S_{2n + 1}, S_{2m} \rangle = \delta_{m, n},\\
    \langle S_{2m}, S_{2n} \rangle = \langle S_{2m + 1}, S_{2n + 1} \rangle = 0.
\end{align*}

Starting this iteration with $S_0$ being a constant function, we retrieve exactly the first $K$ skew-orthogonal polynomials with respect to the inner product $\langle \cdot, \cdot \rangle$ as $S$. $H$ contains the recurrence coefficient as in \eqref{eq:skew_recurrence}.

If we have a finite discrete domain $(x_1, \dots, x_N)$ for our polynomials, we can choose $(f(x_1), \dots, f(x_N))$ as our representation for $f \in V$ and $S$ becomes an $N \times M$ matrix. This means we can rewrite the Arnoldi factorization as
\begin{equation*}
    A S = S H + r \cdot e_{M - 1}^T,
\end{equation*}
where $e_k$ is the $k^\text{th}$ unit vector.

\subsection{Symplectic Gram-Schmidt}

A key step of the proposed symplectic Arnoldi iteration is the skew-orthogonalization procedure to produce a symplectic basis $S$ spanning our Krylov subspace $\{x_0, A x_0, \dots, A^{M - 1} x_0\}$. Matsuo and Nodera \cite{matsuo_block_2014} nicely summarize different symplectic Gram-Schmidt methods.

The basic idea can be illustrated by studying two vectors $x_1$ and $x_2$ one wants to skew-orthogonalize, i.e. factor $X_1 = [x_1, x_2]$ into two matrices $S_1 = [s_1, s_2]$ and $R_1$, such that
\begin{equation*}
    X_1 = S_1 R_1, \quad R_1 = \begin{bmatrix}
        r_{11} & r_{12}\\
        0 & r_{22}
    \end{bmatrix}.
\end{equation*}

With symplectic Gram-Schmidt, one has two degrees of freedom for each pair of vectors, so $r_{11}$ and $r_{12}$ can be chosen arbitrarily and one then computes $s_1$ and $s_2$ as:
\begin{align*}
    s_1 &= x_1 / r_{11} \\
    y &= x_2 - r_{12} s_1 \\
    r_{22} &= \langle s_1, y \rangle \\
    s_2 &= y / r_{22}.
\end{align*}

Matsuo and Nodera describe three common choices for the elementary SR factorization step:
\begin{itemize}
    \item \textsc{esr1}: $r_{11} = \norm{x_1},\ r_{12} = 0$,
    \item \textsc{esr2}: $r_{11} = \norm{x_1},\ r_{12} = s_1^T x_2$,
    \item \textsc{esr3}: $r_{11} = \norm{\langle x_1, x_2 \rangle},\ r_{12} = 0,\ (r_{22} = \pm 1)$.
\end{itemize}

For ease of implementation we introduce a modified \textsc{esr3}, defined as
\begin{itemize}
    \item \textsc{esr3}\small{m}: $r_{11} = 1,\ r_{12} = 0,\ (r_{22} = \langle x_1, x_2 \rangle)$.
\end{itemize}

As opposed to Lanczos, where one just needs to orthogonalize against the previous two vectors, for symplectic Arnoldi, one has to skew-orthogonalize against all previous vectors. The full skew-orthonormalization step is described in \Cref{alg:msgs}.

\begin{algorithm}[t]
\begin{algorithmic}[1]
\Function{skeworthonormalize}{$[s_1, s_2, \dots, s_n], v$}
    \For{$k = 1, \dots, \lfloor n / 2 \rfloor$}
        \State $h_{2k - 1} \gets -\langle s_{2k}, v \rangle$
        \State $h_{2k} \gets \langle s_{2k - 1}, v \rangle$
        \State $v \gets v - h_{2k - 1} s_{2k - 1} - h_{2k} s_{2k}$ \Comment{For CSGS, only update $v \gets v - \sum_{i = 1}^n h_i s_i$ once}
    \EndFor
    \If{$n$ is odd}
        \State $h_n \gets r_{12}(s_n, v)$
        \State $v \gets v - h_n s_n$
        \State $h_{n + 1} \gets \langle s_n, v \rangle$
    \Else
        \State $h_{n + 1} \gets r_{11}(v)$
    \EndIf
    \State $v \gets v / h_{n + 1}$
    \State \Return $v, h$
\EndFunction
\end{algorithmic}
\caption{Modified symplectic Gram-Schmidt step}
\label{alg:msgs}
\end{algorithm}

Here, $r_{11}(x_1)$ and $r_{12}(x_1, x_2)$ depend on the choice of ESR as described above. $h$ is the vector that eventually becomes the $n^\text{th}$ column of the upper Hessenberg matrix $H$ in the Arnoldi factorization.
For classical symplectic Gram-Schmidt (CSGS), $v$ is only updated once in the end, instead of once for every iteration.

\subsection{Reorthogonalization}\label{subsec:reorthogonalization}

In the case of the QR factorization, it has long been known that roundoff errors can accumulate and cause the basis to not have the desired orthogonality properties anymore after a certain number of iterations. For this reason, reorthogonalization has been introduced \cite{daniel_1976}. The same problem also occurs for the symplectic Gram-Schmidt procedure. The full algorithm with iterative renormalization is described in \Cref{alg:reorth}. It does not differ much from the previously described modified Gram-Schmidt algorithm, but instead of iterating through the basis just once, the procedure is repeated until $\norm{v'} \ge \eta \norm{v}$.

In practice, the procedure typically does not need to be run more than twice, so one alternative is to skip the norm check and simply run the procedure twice. This can be slightly more expensive for the first couple of steps, but accuracy and performance for a larger number of steps is very similar.

\begin{algorithm}[h]
\begin{algorithmic}[1]
\Function{skeworthonormalize'}{$[s_1, s_2, \dots, s_n], v; \eta$}
    \State $h \gets 0$
    \State $n_\text{old} \gets \infty$
    \While{$\norm{v} < \eta \cdot n_\text{old}$}
        \State $n_\text{old} \gets \norm{v}$
        \For{$k = 1, \dots, \lfloor n / 2 \rfloor$}
            \State $h_{2k - 1} \gets h_{2k - 1} - \langle s_{2k}, v \rangle$
            \State $h_{2k} \gets h_{2k} + \langle s_{2k - 1}, v \rangle$
            \State $v \gets v - h_{2k - 1} s_{2k - 1} - h_{2k} s_{2k}$
        \EndFor
        \If{$n$ is odd}
            \State $h_n \gets h_n + r_{12}(s_n, v)$
            \State $v \gets v - h_n s_n$
        \EndIf
        \State \Comment{For CSGS, update $v \gets v - \sum_{i = 1}^n h_i s_i$ once here, instead of once for every $k$}
    \EndWhile
    \If{$n$ is odd}
        \State $h_{n + 1} \gets \langle s_n, v \rangle$
    \Else
        \State $h_{n + 1} \gets r_{11}(v)$
    \EndIf
    \State $v \gets v / h_{n + 1}$
    \State \Return $v, h$
\EndFunction
\end{algorithmic}
\caption{Modified symplectic Gram-Schmidt with iterated reorthogonalization}
\label{alg:reorth}
\end{algorithm}

\subsection{The Algorithm}

\begin{algorithm}[h]
\begin{algorithmic}[1]
\Function{arnoldi}{$\langle \cdot , \cdot \rangle: V \times V \rightarrow \mathbb{R}, n$}
    \State $v_0 \gets x^0 / r_{11}(x^0), \quad \Big(\text{e.g. }\underbrace{(r_{11}^{-1}(1), \dots, r_{11}^{-1}(1))}_\text{$m$ times} \text{ for } V = \mathbb{R}^m \Big)$
    \State $S \gets [v_0]$
    \State $H \gets [\,]$
    \For{$k = 1, \dots, n$}
        \State $v_k \gets x \cdot v_{k - 1}$
        \State $(v_k, h) \gets \textsc{skeworthonormalize}(S, v_k)$
        \State $S \gets [S \, |\, v_k]$
        \State $H \gets \begin{bNiceArray}{c|c}
            H & \Block{2-1}{h} \\
            0 \cdots 0
            \CodeAfter
            \tikz \draw (2-|1.1) -- (2-|2);
        \end{bNiceArray}$
    \EndFor
    \State $r \gets x \cdot v_n$
    \State \Return $S, H, r$
\EndFunction
\end{algorithmic}
\caption{The Arnoldi iteration for generating SOPs}
\label{alg:arnoldi}
\end{algorithm}

\Cref{alg:arnoldi} was implemented in \textsc{Julia}~\cite{bezanson_julia_2017} as a skew-orthogonalization routine on top of \texttt{KrylovKit.jl}~\cite{Haegeman_KrylovKit_2024}, extending the existing Arnoldi method to support symplectic Gram-Schmidt. We are working on upstreaming this work to be part of the \texttt{KrylovKit.jl} package. 

\subsection{Numerical stability of symplectic Gram-Schmidt}

We compare our symplectic Arnoldi iteration against computing polynomial coefficients directly using the skew-symmetric variant of the Cholesky decomposition described in \Cref{lem:bjb}. For the Cholesky method, one first builds the matrix
\begin{equation*}
    M = \big( \langle x^i, x^j \rangle_w \big)_{i, j \in \{0, \dots, n - 1\}},
\end{equation*}
which one then factors as $M = B^T J_n B$. The inverse of the upper-triangular matrix $B$ then contains the coefficients of the skew-orthogonal polynomials with respect to $\langle \cdot, \cdot \rangle_w$. One can then use the Vandermonde matrix $V$ to compute the skew-orthogonal basis $W$ as
\begin{equation*}
    W = V(x_0, \dots, x_{n - 1}) \cdot B^{-1},
\end{equation*}
where the $i^\text{th}$ column of $W$ is the $i^\text{th}$ basis element evaluated at each $x_j$.

\begin{figure}
    \centering
    \includegraphics[width=0.7\textwidth]{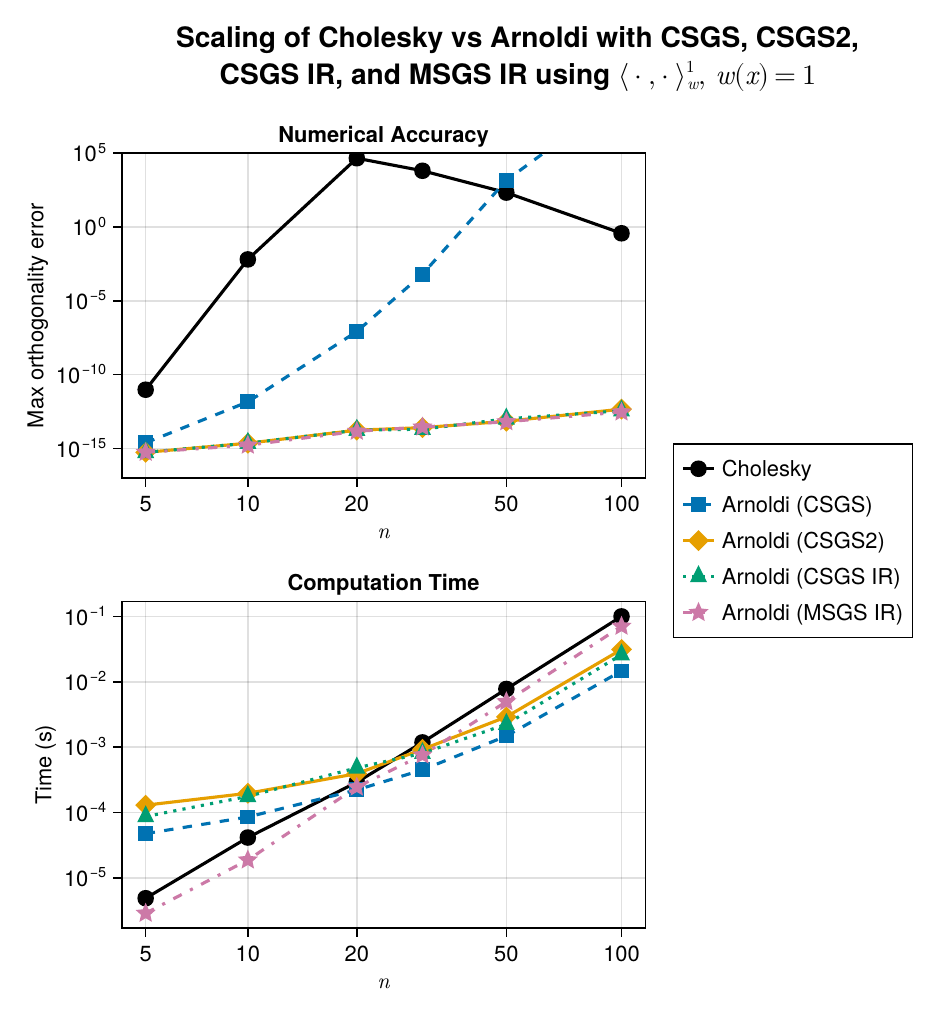}
    \caption{ 
    Comparing different methods for computing skew-orthogonal polynomial bases, it is evident that our symplectic Arnoldi method significantly outperforms the method based on the skew-symmetric Cholesky factorization (in black) in terms of numerical accuracy. Classical symplectic Gram-Schmidt without reorthogonalization (CSGS) also failed to maintain skew-orthogonality for larger $n$, while
    classical symplectic Gram-Schmidt, reorthogonalizing once at every iteration (CSGS2), classical symplectic Gram-Schmidt using iterated reorthogonalization (CSGS IR) and modified symplectic Gram-Schmidt, also using iterated reorthogonalization (MSGS IR) all delivered similar numerical accuracy. MSGS IR was more efficient for small $n$, while CSGS2 and CSGS IR performed better for large $n$ due to their ability to be parallelized. There was not much difference between CSGS2 and CSGS IR, though CSGS IR had a slight performance advantage over CSGS2 in some tests.
    In all cases, ESR2 was chosen as the normalization and for the iterated reorthogonalization methods, $\eta$ was chosen to be 0.75.}
    \label{fig:scaling_cholesky_arnoldi_one}
\end{figure}

As we can see in \Cref{fig:scaling_cholesky_arnoldi_one}, for $\beta = 1$ with uniform weight ($w(x) = 1$), this method is much less numerically stable than our symplectic Arnoldi iteration. One can also see the necessity of the reorthogonalization procedure described in \Cref{subsec:reorthogonalization}, as classical symplectic Gram-Schmidt without reorthogonalization (CSGS) quickly diverges.

All the other methods with reorthogonalization performed similarly accuracy-wise, though there were performance differences. For small $n$, modified Gram-Schmidt (MSGS) was generally the fastest, while for large $n$, CSGS can take advantage of parallelization and is generally faster than MSGS. The difference between reorthogonalizing once at every iteration (CSGS2) and iterated reorthogonalization (CSGS IR) were almost negligible, with CSGS2 perhaps having a slight edge over CSGS IR for very small $n$.

\begin{figure}
    \centering
    \includegraphics[width=0.7\textwidth]{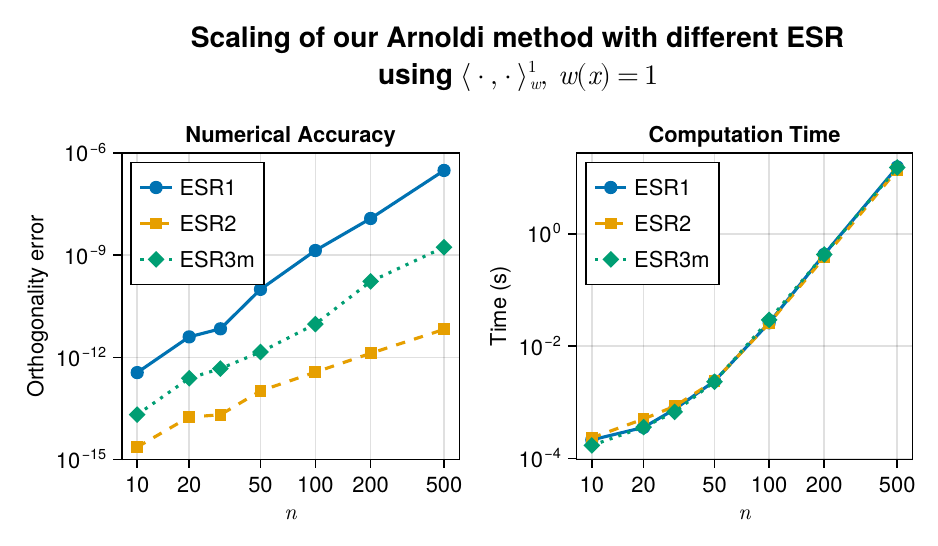}
    \caption{Comparison of different normalizations for skew-orthogonalization. In all cases, classical symplectic Gram-Schmidt with iterated reorthogonalization (CSGS IR) and $\eta = 0.75$ was used}
    \label{fig:scaling_esr_one}
\end{figure}

\begin{figure}
    \centering
    \includegraphics[width=0.7\textwidth]{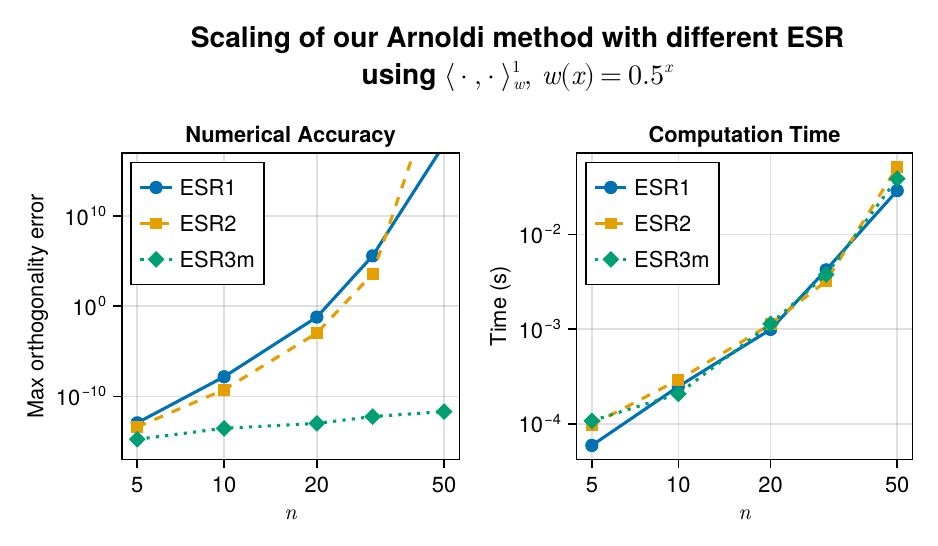}
    \caption{Comparison of different normalizations for skew-orthogonalization on an exponentially decaying weight function. In all cases, classical symplectic Gram-Schmidt with iterated reorthogonalization (CSGS IR) and $\eta = 0.75$ was used}
    \label{fig:scaling_esr_power}
\end{figure}

Again, for $\beta = 1$ with uniform weight ($w(x) = 1$), ESR2 is the most stable (see \Cref{fig:scaling_esr_one}), as already mentioned for the typical case by Matsuo and Nodera \cite{matsuo_block_2014}. For many PfPPs that we are interested in, the weight decays exponentially or faster, so we also ran tests with $w(x) = 0.5^x$, see \Cref{fig:scaling_esr_power}. There, ESR3m vastly outperformed ESR1 and ESR2, which both have unacceptable max and mean errors beyond $n = 20$. This is most likely due to the $L^2$ norm of all the even vectors exploding, causing numerical instabilities.

All benchmarks were run on a 
laptop with a \textsc{Ryzen AI 7 Pro 360} processor, containing 3 Zen 5 and 5 Zen 5c cores, running at a maximum of 5 and 3.3\,GHz, respectively.
Plots were created using the \texttt{Makie.jl}~\cite{DanischKrumbiegel2021} plotting library.

\section{Numerical Experiments}\label{sec:numerical}

In this section, we discuss various applications of Pfaffian point processes and numerical experiments related to them. 

\subsection{Symmetric corner growth}

Johannson \cite{johansson_shape_2000} introduced a symmetric corner growth model analogously to the $\beta = 2$ standard geometric growth. The restriction of the standard corner growth model to the symmetric one is analogous to the restriction of the permutations to involutions. 

\begin{figure}[h]
    \centering
    \includegraphics[width=0.6\textwidth]{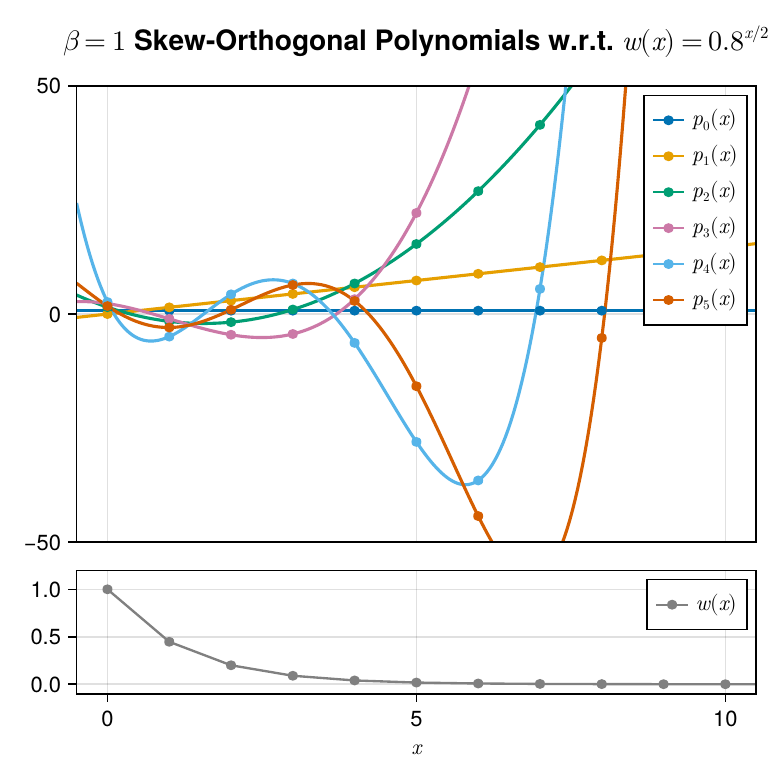}
    \caption{Shown above are the skew-orthogonal polynomials with respect to the discrete weight $w(x) = 0.8^{x / 2}, \ x \in \mathbb{Z}_{\ge 0}$, corresponding to symmetric corner growth with $p = 0.8$. Below is the weight function itself.}
    \label{fig:symmetric_growth_sop}
\end{figure}

Let $A = \big(w(i, j)\big)$ be a waiting time matrix with random entries according to the following distribution:
\begin{equation*}
    w(i, j) \sim \begin{cases}
        \operatorname{Geometric}(1 - \sqrt{q}) & \text{if } i = j \\
        \operatorname{Geometric}(1 - q) & \text{if } i < j \\
        = w(j, i) & \text{if } i > j
    \end{cases}.
\end{equation*}

Define $\Pi_{N, N}$ as the set of all up-right paths from $(1, 1)$ to $(N, N)$, then $F(N)$ is defined as:
\begin{equation*}
    F(N) = \max \Big\{ \sum_{(i, j) \in \pi} w(i, j) : \pi \in \Pi_{N, N} \Big\}.
\end{equation*}

Johannson derives its distribution as
\begin{equation*}
    \mathbb{P}[F(N) \le t] = \frac{1}{Z_N^{(1)}} \sum_{\substack{h \in \mathbb{N}^N \\ \max \{h_j\} \le t + N - 1}} \prod_{1 \le i < j \le N} |h_j - h_i| \prod_{i = 1}^N q^{h_i / 2}.
\end{equation*}

This is exactly of the form of the discrete Coulomb gas for $\beta = 1$, so we can construct $K$ from skew-orthogonal polynomials with the $\beta = 1$ skew-inner product $\langle \cdot, \cdot \rangle_w^1$, taking $w(x) = q^{x / 2}$ as weight function. These are shown in \Cref{fig:symmetric_growth_sop}.

\begin{figure}
    \centering
    \includegraphics[width=0.6\textwidth]{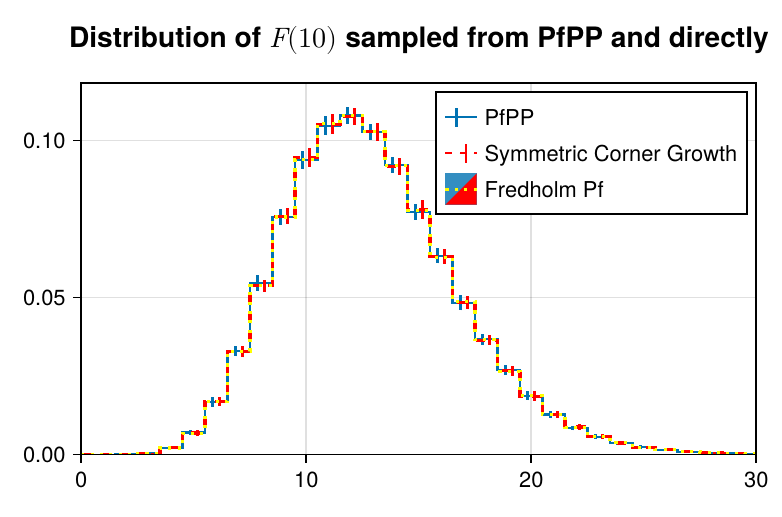}
    \caption{The distribution of $F(10)$, the time it takes for symmetric corner growth to reach the square at $(10, 10)$, alongside the same distribution sampled from the corresponding PfPP.}
    \label{fig:symmetric_growth_pfpp}
\end{figure}

For sampling the DPP we have to choose a cutoff, ideally at a point where the weight function is almost zero. We can also cut it off earlier though, as long as we reject any $F(N) > n_\text{cutoff} - M + 1$. The sampled distribution of $F(10)$ can be seen in \Cref{fig:symmetric_growth_pfpp}.

\subsection{Finite $N$ Random Matrices}

We consider random matrices with finite size. In particular, we consider the eigenvalues of GOE and GSE. Obviously, the GOE and GSE can be easily sampled by using their matrix models, and for large $N$, one can also use the tridiagonal model \cite{dumitriu2002matrix}. 

\subsubsection{Finite $N$ GOE}\label{sec:GOEfinite}

We start with an $N \times N$ matrix $X_N$ with real i.i.d entries sampled from a standard normal distribution. Symmetrize the matrix and normalize the result to obtain the matrix $W_N^{(1)}$ as follows
\begin{equation*}
    X_{N, ij} \sim \mathcal{N}(0, 1), \quad
    W_N^{(1)} = \frac{1}{\sqrt{2}} \left( X_N + X_N^T \right).
\end{equation*}
The eigenvalues of $W_N^{(1)}$ have the joint density
\begin{equation*}
    f_N^{(1)}(\lambda_1, \dots, \lambda_N) = \frac{1}{\mathcal{Z}_{N, 1}} \prod_{1 \le i < j \le N} |\lambda_j - \lambda_i| \prod_{i = 1}^N e^{-\frac{1}{4} \lambda_i^2},
\end{equation*}
where $\mathcal{Z}_{N, 1}$ can be computed from \eqref{eq:partitionfunction}. 

We now want to find the skew-orthogonal polynomials with respect to the inner product \eqref{eq:beta1innerproduct}
\begin{equation*}
    \langle f, g \rangle_\text{GOE} \coloneq \frac{1}{2} \int_{-\infty}^\infty \int_{-\infty}^\infty f(x) g(y) \cdot \operatorname{sign}(y - x) e^{-\frac{1}{4} x^2} e^{-\frac{1}{4} y^2} \dd{x} \dd{y}.
\end{equation*}
These skew-orthogonal polynomials have already been derived analytically, see \cite{mehta2004random}, but we recompute them numerically in order to demonstrate our symplectic Arnoldi method on continuous kernels.

Since adaptive quadrature methods struggle with the discontinuity at $x = y$, the inner integral is split in two and both triangular integrals are evaluated separately:
\begin{equation*}
    \langle f, g \rangle_\text{GOE} = \frac{1}{2} \int_{-\infty}^\infty f(x) e^{-\frac{1}{4} x^2} \left( \int_{x}^\infty g(y) e^{-\frac{1}{4} y^2} \dd{y} - \int_{-\infty}^x g(y) e^{-\frac{1}{4} y^2} \dd{y} \right) \dd{x}
\end{equation*}

\begin{figure}
    \centering
    \includegraphics[width=0.6\textwidth]{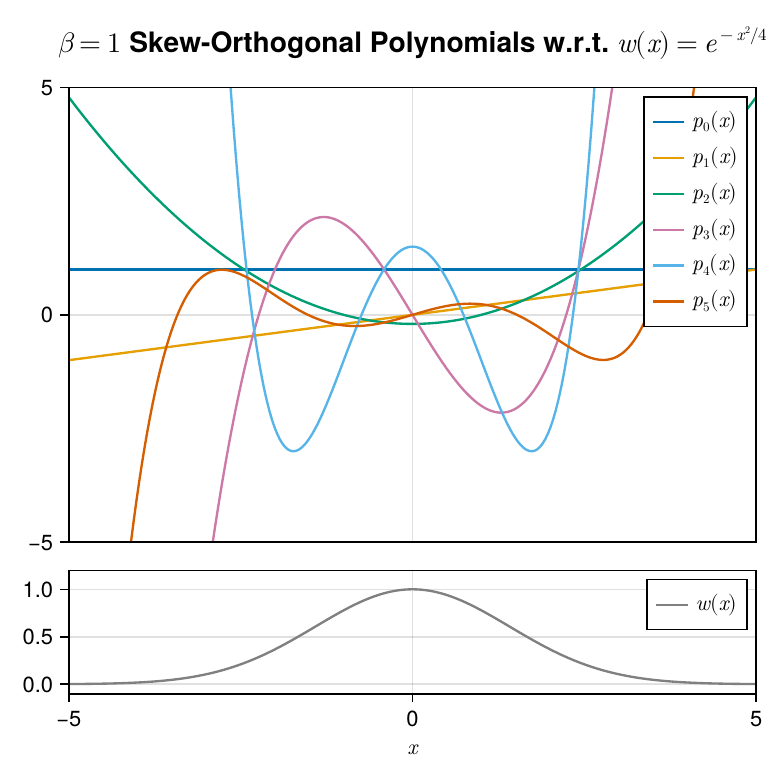}
    \caption{The first six skew-orthogonal polynomials with respect to the GOE skew-inner product.}
    \label{fig:goe_sop}
\end{figure}

To represent the polynomials without having to discretize the domain, \texttt{ApproxFun.jl}~\cite{ApproxFun.jl-2014} was used with the standard ($\beta = 2$) Hermite polynomials chosen as a basis.
The integration for the skew-inner product was performed using \texttt{QuadGK.jl}~\cite{quadgk}, interfacing through \texttt{Integrals.jl}~\cite{DifferentialEquations.jl-2017} and \texttt{KrylovKit.jl}~\cite{Haegeman_KrylovKit_2024}, extended with our symplectic Gram-Schmidt procedure, was used for the Arnoldi iteration.
This results in $N$ skew-orthogonal polynomials as the symplectic basis of the Kryvlov subspace $\{1, x, x^2, \dots, x^{N - 1}\}$ with respect to the skew-inner product $\langle \cdot, \cdot \rangle_\text{GOE}$, the first 6 of which are shown in \Cref{fig:goe_sop}.

As derived by Metha in \cite{mehta2004random}, in the case of even $N$, the Pfaffian kernel $K_N$ can be defined based on these SOPs $R_k(x)$ as follows:

Define
\begin{align*}
    w(x) &= e^{-\frac{1}{4} x^2} \\
    \psi_k(x) &= \frac{1}{2} \int_{-\infty}^\infty R_k(y) \cdot \operatorname{sign}(x - y) w(y) \dd{y},
\end{align*}
as well as
\begin{align*}
    S_1(x, y) &= w(x) \cdot \sum_{k = 0}^{N / 2 - 1} \big( R_{2k + 1}(x) \psi_{2k}(y) - R_{2k}(x) \psi_{2k + 1}(y) \big), \\
    D_1(x, y) &= w(x) w(y) \cdot \sum_{k = 0}^{N / 2 - 1} \big( -R_{2k + 1}(x) R_{2k}(y) + R_{2k}(x) R_{2k + 1}(y) \big), \\
    J_1(x, y) &= \sum_{k = 0}^{N / 2 - 1} \big( \psi_{2k + 1}(x) \psi_{2k}(y) - \psi_{2k}(x) \psi_{2k + 1}(y) \big) - \frac{1}{2} \operatorname{sign}(x - y).
\end{align*}
Then, the $2\times 2$ matrix valued kernel for the PfPP is defined as
\begin{equation*}
    K^{(1)}_N(x, y) = \begin{bNiceArray}{rr}
        J_1(x, y) & S_1(y, x) \\ -S_1(x, y) & -D_1(x, y)
    \end{bNiceArray}.
\end{equation*}

\begin{figure}[h]
    \centering
    \includegraphics[width=0.6\textwidth]{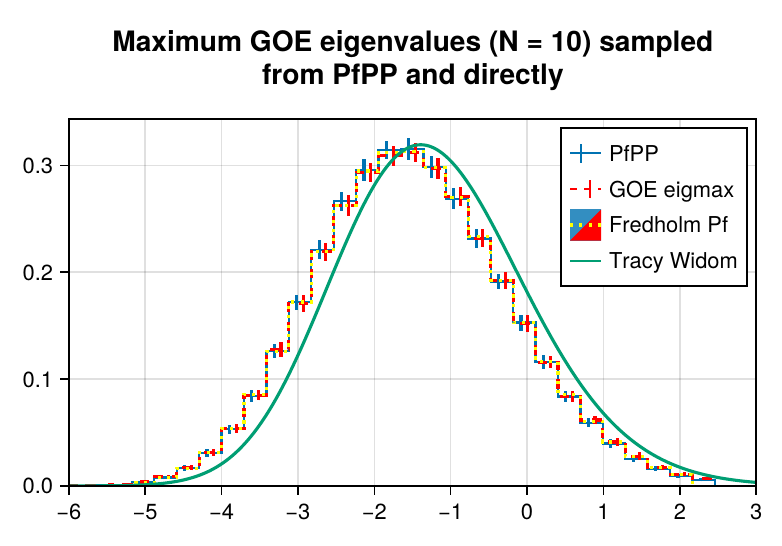}
    \caption{Shown above is the distribution of the maximum eigenvalues of the $N = 10$ GOE, normalized as described in \eqref{eq:softedgescaling} to match the $\beta = 1$ Tracy-Widom distribution as $N \rightarrow \infty$. The solid blue lines and error bars correspond to the values sampled from the discretized PfPP using our modified skew-Cholesky method. The dashed red line corresponds to a histogram of the maximum eigenvalues directly sampled from the GOE matrices. The dotted yellow line is the Fredholm Pfaffian calculated from our discretized PfPP and the solid green line is the $\beta = 1$ Tracy-Widom distribution. For both the samples from the PfPP as well as the direct samples, 50 batches of 10000 samples were used. These batches were then averaged for each bin, with the error bars representing the standard deviation between batches.}
    \label{fig:goe_hist}
\end{figure}

Note how $K^{(1)}_N(x, y)$ deviates from the $K_1(x, y)$ as defined by Mehta. This is due to Mehta defining the PfPP in terms of the quaternion determinant $\operatorname{qdet} K_1(X, X) \coloneq \pf{J \cdot K_1(X, X)}$, where as we define the PfPP directly in terms of the Pfaffian of $K^{(1)}_N(X, X)$. We also work with normalized SOPs instead of monic ones, so $q_k = 1$.

We sample from this kernel using \Cref{alg:pfppsampler} and compare against GOE eigenvalues sampled directly from symmetric random matrices. Histograms for the largest eigenvalue can be seen in \Cref{fig:goe_hist}, alongside the Fredholm Pfaffian as well as the $\beta = 1$ Tracy-Widom PDF, which this distribution will converge to for $N \rightarrow \infty$.

\begin{remark}[Numerical difficulties for orthogonal ensembles]
For any given $\beta=1$ polynomial ensemble, we encounter a discontinuity on the diagonal from $\epsilon(x,y)$ term. When one is computing the Fredholm Pfaffians, the discontinuity induces numerical inaccuracy of the discretization. Similarly, when we discretize the continuous PfPP to a discrete one, we lose accuracy. When computing the Fredholm Pfaffian, Bornemann \cite{bornemann2010numerical} discusses an alternative method using the relationship between orthogonal and unitary ensembles to remedy this problem. 
\end{remark}

\subsubsection{Finite $N$ GSE Eigenvalues}\label{sec:GSEfinite}

For the Gaussian symplectic ensemble (GSE), start with two $N \times N$ matrices $X_N$ and $Y_N$ with complex i.i.d entries sampled from a standard normal distribution. Obtain the matrix $W_N^{(4)}$ via:
\setlength{\extrarowheight}{1mm}
\begin{align*}
    (X_N)_{ij}, (Y_N)_{ij} \sim \mathcal{N}(0, 1 / 2) + i \mathcal{N}(0, 1 / 2), \quad
    A_N = \begin{bNiceArray}{rc}
        X & Y \\
        -\overline{Y} & \overline{X}
    \end{bNiceArray}, \quad
    W_N^{(4)} = \frac{1}{\sqrt{8}} \left( A_N + A_N^* \right).
\end{align*}
\setlength{\extrarowheight}{0mm}
Alternatively, one can take a quaternion $N\times N$ matrix $A_N$ with i.i.d entries from
\begin{equation*}
     \mathcal{N}(0, 1 / 2) + i \mathcal{N}(0, 1 / 2) +j \mathcal{N}(0, 1 / 2) + k \mathcal{N}(0, 1 / 2),
\end{equation*}
where $i, j, k$ are the imaginary units of the quaternions. Then, one can take the matrix $(A_N + A_N^*)/\sqrt{8}$ as the GSE. Since this matrix $W_N^{(4)}$ is Hermitian (or self-dual if using the quaternionic version), the eigenvalues are real, occur in pairs and have a joint density of:
\begin{equation*}
    f_N^{(4)}(\lambda_1, \dots, \lambda_N) = \frac{1}{\mathcal{Z}_{N, 4}} \prod_{1 \le i < j \le N} |\lambda_j - \lambda_i|^4 \prod_{i = 1}^N e^{-\lambda_i^2}.
\end{equation*}

\begin{figure}[b]
    \centering
    \includegraphics[width=0.6\textwidth]{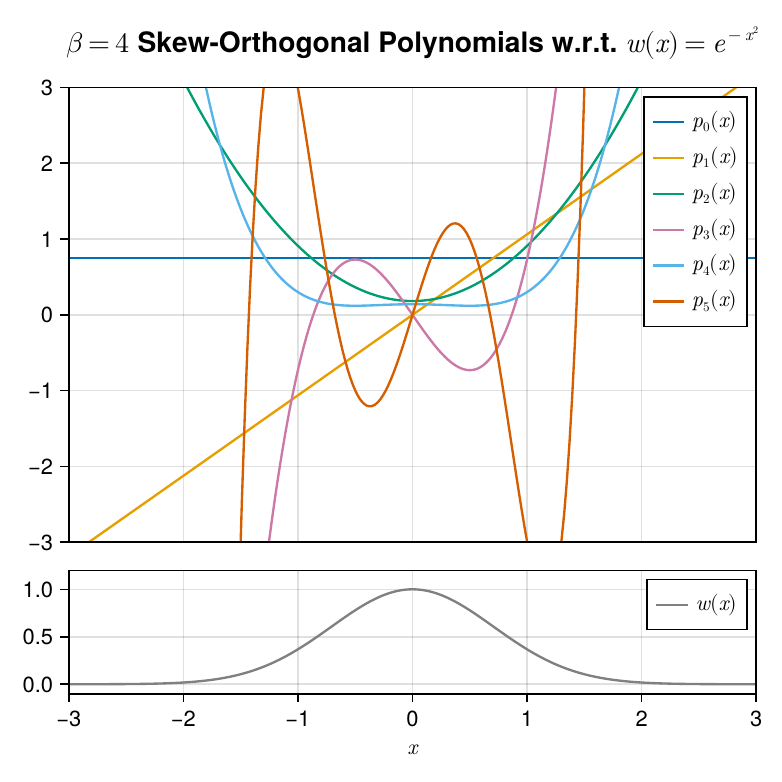}
    \caption{The first six skew-orthogonal polynomials with respect to the GSE skew-inner product \eqref{eq:gseinnerproduct}.}
    \label{fig:gse_sop}
\end{figure}

As in the previous examples we want to find the corresponding skew-orthogonal polynomials, this time with respect to the skew-inner product \eqref{eq:beta4innerproduct}
\begin{equation}\label{eq:gseinnerproduct}
    \langle f, g \rangle_\text{GSE} \coloneq \int_{-\infty}^\infty \left( f(x) g'(x) - f'(x) g(x) \right) \cdot e^{-2x^2} \dd{x}.
\end{equation}
We proceed similarly to GOE for constructing these polynomials. The resulting polynomials can be seen in \Cref{fig:gse_sop}.

Define
\begin{align*}
    w(x) &= e^{-x^2} \\
    S_4(x, y) &= w(x) w(y) \cdot \sum_{k = 0}^{N - 1} \big( Q'_{2k + 1}(x) Q_{2k}(y) - Q'_{2k}(x) Q_{2k + 1}(y) \big), \\
    D_4(x, y) &= w(x) w(y) \cdot \sum_{k = 0}^{N - 1} \big( -Q'_{2k + 1}(x) Q'_{2k}(y) + Q'_{2k}(x) Q'_{2k + 1}(y) \big), \\
    I_4(x, y) &= w(x) w(y) \cdot \sum_{k = 0}^{N - 1} \big( Q_{2k + 1}(x) Q_{2k}(y) - Q_{2k}(x) Q_{2k + 1}(y) \big),
\end{align*}
then the PfPP kernel is defined as:
\begin{equation*}
    K^{(4)}_N(x, y) = \begin{bNiceArray}{rr}
        I_4(x, y) & S_4(y, x) \\ -S_4(x, y) & -D_4(x, y)
    \end{bNiceArray}.
\end{equation*}

\begin{figure}[h]
    \centering
    \includegraphics[width=0.6\textwidth]{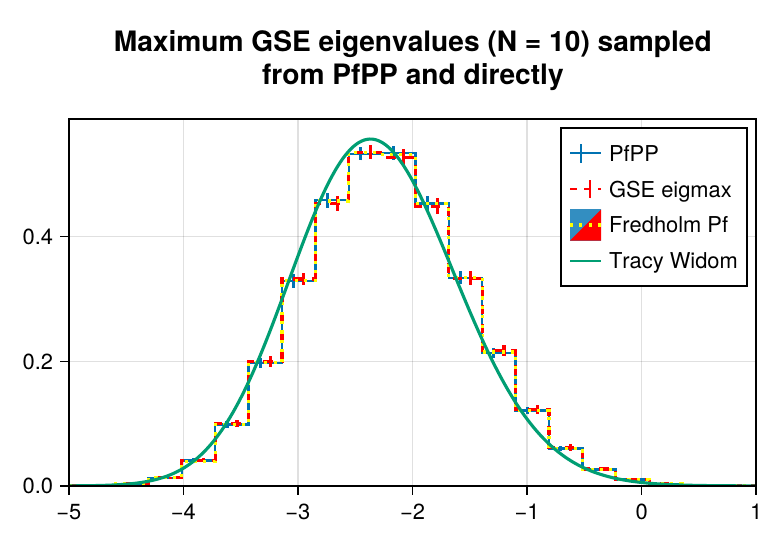}
    \caption{The above histograms were again produced as already described for \Cref{fig:goe_hist}, by simply replacing the PfPP with the one calculated from the GSE skew-inner product and the maximum eigenvalues taken from GSE instead of GOE matrices. As discussed there, 50 batches of 10000 samples were used for calculating the mean and standard deviation of each bin. The normalization is $\lambda \mapsto (2N)^{1/6}(\sqrt{2}\lambda - 2\sqrt{N})$, to match the convention of the Tracy-Widom distribution.}
    \label{fig:gse_hist}
\end{figure}

Besides the differences from Mehta's definition already mentioned in the case of the GOE, we also work with $w(x) = e^{-x^2}$ instead of $w(x) = e^{-2x^2}$, therefore the missing square root over $w(x) w(y)$.

We proceed, as we did for GOE, with sampling from the PfPP using our skew-Cholesky method and compare against the GSE eigenvalues. The results can be seen in \Cref{fig:gse_hist}.

\subsubsection{Continuous Sampling}

\begin{figure}[b]
    \centering
    \includegraphics[width=0.6\textwidth]{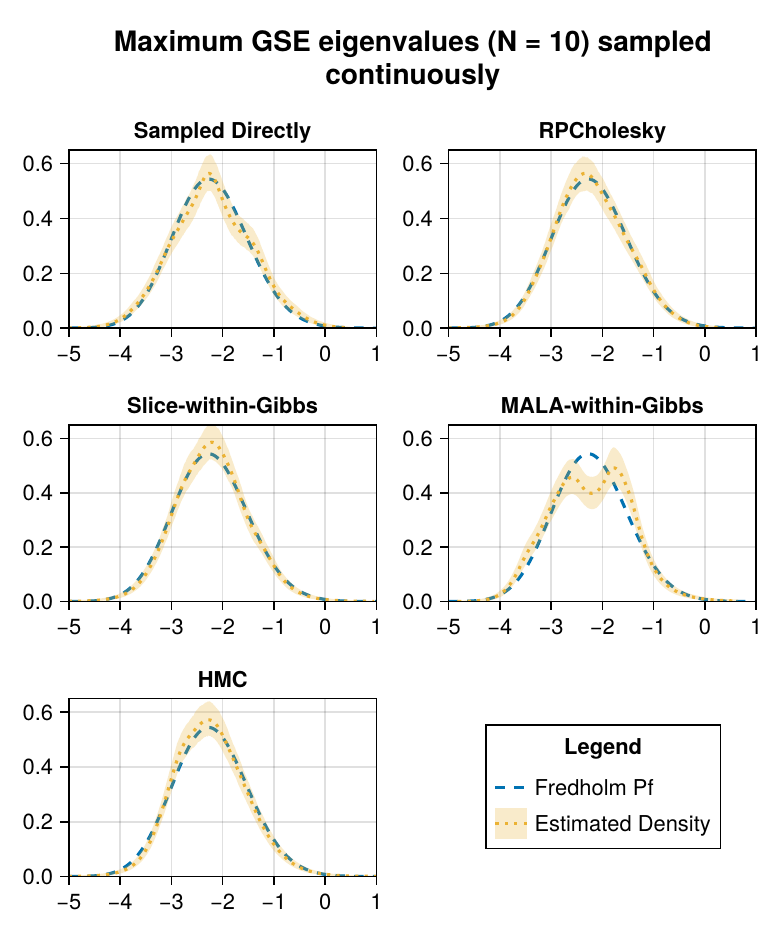}
    \caption{Comparison of different continuous sampling algorithms from 1000 samples taken for each of the methods. All methods, except for \textsc{MALA-within-Gibbs}, show similarly accurate results, which struggles due to strong autocorrelation between samples due to the nature of the Metropolis-adjusted Langevin algorithm. This means \textsc{MALA-within-Gibbs} requires a larger number of samples to give accurate density estimates. The package \texttt{BayesDensity.jl} \cite{simensen2026bayesdensity} was used for density estimation using variational inference to estimate the posterior of a \texttt{HistSmoother}, except for the \textsc{MALA-within-Gibbs} samples, where augmented Gibbs sampling had to be used due to convergence issues. Normalization is the same as in \Cref{fig:gse_hist}.}
    \label{fig:gse_continuous}
\end{figure}

For the GSE, we also tested the continuous sampling methods as described in \Cref{sec:continuous}, see \Cref{fig:gse_continuous}.
1000 samples were produced using each of the methods and from those samples, the density was estimated using \texttt{BayesDensity.jl} \cite{simensen2026bayesdensity}. For the two Gibbs sampling methods, 100 additional burn-in samples were used to avoid bias due to the initial guess. HMC used \texttt{DynamicHMC.jl}'s \cite{tamas_k_papp_2026_18130162} default parameters, which already includes burn-in. One can see that \textsc{MALA-within-Gibbs} struggles with convergence issues, due to the fact that samples don't necessarily get updated on every step, so there is a strong autocorrelation between samples. We verified convergence does improve with a higher number of samples, but not enough to be compensated by the slight performance advantage over \textsc{Slice-within-Gibbs}.

\subsection{Airy point processes}

In this section we use \Cref{alg:pfppsampler} to compute statistics of $\beta=1,4$ Airy point processes. We first briefly review the $\beta=2$ case. For an $N\times N$ GUE defined with \eqref{eq:hermitebeta} with $\beta=2$, recenter and rescale the eigenvalues as $N\to\infty$ with
\begin{equation}\label{eq:softedgescaling}
    \lambda \longmapsto N^{1/6}(\lambda - 2\sqrt{N}).
\end{equation}
This map corresponds to removing the first order term of the largest eigenvalue $\lambda_1$ and properly rescaling. The scaling \eqref{eq:softedgescaling} is often called the soft-edge scaling limit. 

In \cite{tracy1994level}, Tracy and Widom proved that the distribution function $F_2(s)$ of the largest eigenvalue under the above scaling limit can be expressed in terms of the solution to the Painleve II ordinary differential equation. Moreover, it is also known that the point process of (countably many) eigenvalues under this scaling limit is a DPP with the Airy kernel 
\begin{equation*}
    K_{\Ai}(x, y) = \frac{\Ai(x)\Ai'(y) - \Ai'(x)\Ai(y)}{x-y}. 
\end{equation*}
This point process is also called the Airy point process or Airy$_2$ point process to emphasize the $\beta$ parameter. 

The same scaling limit of \eqref{eq:softedgescaling} can be considered for the Hermite (Gaussian) ensemble with any $\beta>0$. The scaled eigenvalues are sometimes denoted by Airy$_\beta$ point process and the distribution of the largest particle is often denoted by $F_\beta(s)$.

\subsubsection{GOE scaling limit}
\begin{figure}[b]
    \centering
    \includegraphics[width=0.6\textwidth]{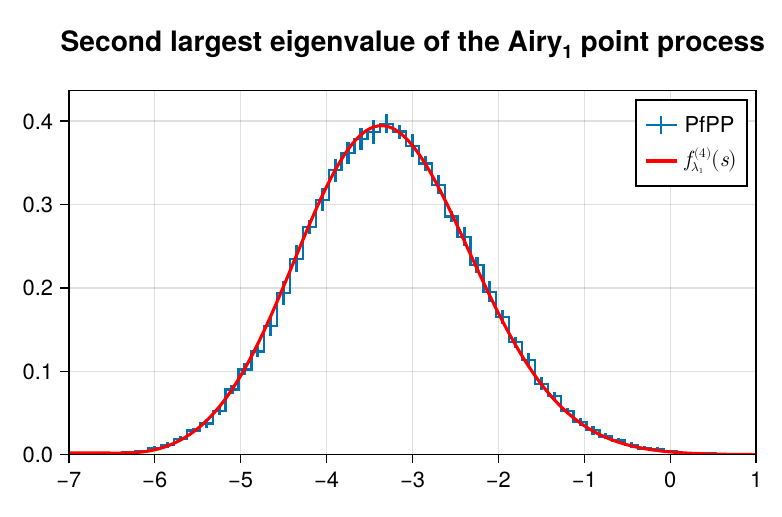}
    \caption{The density of the second largest eigenvalue of the Airy$_1$ point process. The histogram is obtained from \Cref{alg:pfppsampler}, and the red line is obtained from the formula for
    $f_{\lambda_2}^{(1)}(s)$ described in~\eqref{eq:secondeigpdf}.}
    \label{fig:goeairy1}
\end{figure}
With a similar derivation to its finite counterpart (\Cref{sec:GOEfinite}), we have the following $2\times 2$ matrix valued kernel $K^{(1)}$ for the Airy$_1$ point process \cite{a2005matrix}:
\begin{align*}
    K_{11}^{(1)}(x, y) &= \int_0^\infty \Ai(x+z)\Ai'(y+z) \dd{z}- \int_0^\infty \Ai'(x+z)\Ai(y+z)\dd{z}, \\
    K_{12}^{(1)}(x, y) &=  \int_0^\infty\Ai(x+z)\Ai(y+z)\dd{z} + \frac{1}{2}\Ai(x)\int_0^\infty \Ai(y-z)\dd{z}, \\
    K_{21}^{(1)}(x, y) &= -K_{12}^{(1)}(y,x), \\
    K_{22}^{(1)}(x, y) &= \frac{1}{4}\int_0^\infty \int_\lambda^\infty \Ai(x-z)\Ai(y-w)\dd{w} \dd{z}.
\end{align*}
Note that the kernel for this PfPP is not unique, but all kernels are equivalent in the sense that they induce the same Pfaffians and correlation function.

We perform \Cref{alg:pfppsampler} to sample the second largest eigenvalue $\lambda_2$ of the Airy$_1$ point process and compare the histogram with the computed probability density of $\lambda_2$ obtained by the conditional kernel approach introduced in \cite{edelman2023conditional}. We modify Equation (3.21) of \cite{edelman2023conditional} to get the following formula for the density of the second eigenvalue of the Airy$_1$ point process (which also works for $\beta=4$)
\begin{equation}\label{eq:secondeigpdf}
    f_{\lambda_2}^{(1)}(s) = K_{12}^{(1)}(s, s) \cdot \operatorname{skewtr}\left(J K^{(1, s)}(J - K^{(1, s)})\right)\!\!\restriction_{(s, \infty)} \cdot \pf(J - K^{(1, s)}),
\end{equation}
where the skew-trace $\operatorname{skewtr}$ is defined as
\begin{equation*}
    \operatorname{skewtr}(K)\!\!\restriction_{\Omega} \ \coloneq \int_\Omega K_{12}(x,x) dx,
\end{equation*}
for the kernel of the integral operator $K$, the conditional kernel $K^{(1,s)}$ is defined as
\begin{equation*}
    K^{(1,s)} = K^{(1)}(x,y) - K^{(1)}(x,s)\left(K^{(1)}(s,s)\right)^{-1}K^{(1)}(s,y),  
\end{equation*}
and $J$ is defined as $2\times 2$ matrix valued kernel with $J(x,y) = J_1 \delta_{xy}$. Note that the kernel $K^{(1,s)}$ is the Pfaffian point process kernel of the Airy$_1$ point process conditioned on having an eigenvalue at $s$, which is studied in \cite{bufetov2021conditional}. The product of the second and third term in \eqref{eq:secondeigpdf} is the probability of only one eigenvalue existing in $(s, \infty)$ given an eigenvalue at $s$. See Figure~\ref{fig:goeairy1} for the comparison.

\subsubsection{GSE scaling limit}

The scaling limit of the GSE at the soft-edge (Airy$_4$ point process) is a PfPP with the following $2\times 2$ matrix valued kernel \cite{a2005matrix,baik2018pfaffian}
\begin{align*}
    K_{12}^{(4)}(x, y) &= \frac{1}{2}K_{\Ai}(x, y) - \frac{1}{4}\Ai(x) \int_y^\infty \Ai (z)\dd{z},\\
    K_{12}^{(4)}(x, y) &= -\frac{1}{2}\frac{\partial}{\partial y} K_{\Ai}(x, y) - \frac{1}{4}\Ai(x)\Ai(y),\\
    K_{21}^{(4)}(x, y) &= -K_{12}^{(4)}(y,x),\\
    K_{22}^{(4)}(x, y) &= -\frac{1}{2}\int_x^\infty K_{\Ai}(z, y)\dd{z} + \frac{1}{4}\int_x^\infty \Ai(z)\dd{z} \int_y^\infty \Ai(w)\dd{w}.
\end{align*}

\begin{figure}[h]
    \centering
    \includegraphics[width=0.6\textwidth]{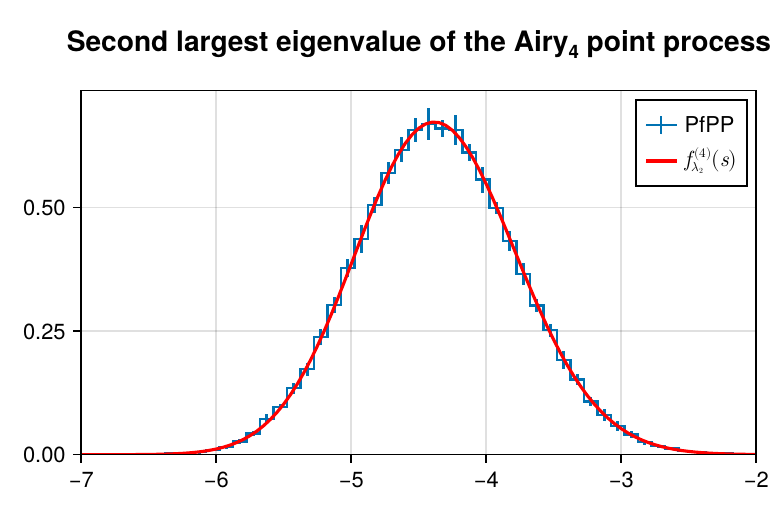}
    \caption{The density of the second largest eigenvalue of the Airy$_4$ point process. The histogram is obtained from \Cref{alg:pfppsampler}, and the red line is obtained from the formula for
    $f_{\lambda_2}^{(4)}(s)$, which is defined in the same as $f_{\lambda_2}^{(1)}(s)$ in~\eqref{eq:secondeigpdf}.}
    \label{fig:gseairy4}
\end{figure}

Once more, we sample the second largest eigenvalue of Airy$_4$ point process using \Cref{alg:pfppsampler} and compare the histogram against its density computed with \eqref{eq:secondeigpdf} in Figure~\ref{fig:gseairy4}.

\subsection*{Acknowledgements}
We thank Jutho Haegeman for his review and feedback on the \textsc{Julia} implementation for. The first author is supported by the U.S. National Science Foundation under award Nos CNS-2346520,  RISE-2425761, and DMS-2325184, by the Defense Advanced Research Projects Agency (DARPA) under Agreement No. HR00112490488, by the Department of Energy, National Nuclear Security Administration under Award Number DE-NA0003965 and by the United States Air Force Research Laboratory under Cooperative Agreement Number FA8750-19-2-1000. Neither the United States Government nor any agency thereof, nor any of their employees, makes any warranty, express or implied, or assumes any legal liability or responsibility for the accuracy, completeness, or usefulness of any information, apparatus, product, or process disclosed, or represents that its use would not infringe privately owned rights. Reference herein to any specific commercial product, process, or service by trade name, trademark, manufacturer, or otherwise does not necessarily constitute or imply its endorsement, recommendation, or favoring by the United States Government or any agency thereof.

\bibliographystyle{alpha}
\bibliography{references}

\end{document}